\documentclass[12pt]{article}

\usepackage{amssymb}

\makeatletter
 
 \@addtoreset{equation}{section}
\makeatother

\newcommand{\ul}{\underline}

\newcommand{\li}{\mbox{\rm li}\,}

\newcommand{\Gal}{\mbox{\rm Gal}}
\newcommand{\Aut}{\mbox{\rm Aut}}
\newcommand{\id}{\mbox{\rm id}}
\newcommand{\Frac}{\displaystyle\frac}
\newcommand{\psum}{\mathop{{\sum}'}}
\newcommand{\scr}{\scriptstyle}
\newcommand{\scrs}{\scriptsize}
\renewcommand{\pmod}[1]{%
\ (\mbox{\rm mod}\ #1)}
\newcommand{\lcm}[1]{%
\langle #1 \rangle}

\newcommand{\prf}{\noindent{\bf Proof. }}

\newcommand{\qed}{\hbox{\rule[-2pt]{5pt}{11pt}}}

\newtheorem{thm}{Theorem}[section]
\newtheorem{prop}{Proposition}[section]
\newtheorem{cor}{Corollary}[section]
\newtheorem{lem}{Lemma}[section]

\begin{document}
\title{On a distribution property of the residual order of $a\pmod p$ --- II}
\author{Leo Murata\footnotemark[1] and Koji Chinen\footnotemark[7]}
\date{}
\maketitle

\begin{center}
{\it Dedicated to P.D.T.A. Elliott on the occasion of his 60th birthday.}
\end{center}

\begin{abstract}
 Let $a$ be a positive integer which is not a perfect $h$-th power with $h\geq2$, and $Q_a(x;4,l)$ be the set of primes $p\leq x$ such that the residual order of $a\pmod p$ in ${\bf Z}/p{\bf Z}^\times$ is congruent to $l$ modulo 4. When $l=0,2$, it is known that calculations of $\sharp Q_a(x;4,l)$ are simple, and we can get their natural densities unconditionally. On the contrary, when $l=1,3$ , the distribution properties of $Q_a(x;4,l)$ are rather complicated. In this paper, which is a sequel of our previous paper \cite{CM}, under the assumption of Generalized Riemann Hypothesis, we determine completely the natural densities of $\sharp Q_a(x;4,l)$ for $l=1,3$.
\end{abstract}
\footnotetext[1]{Department of Mathematics, Faculty of Economics, Meiji Gakuin University, 1-2-37 Shirokanedai, Minato-ku, Tokyo 108-8636, Japan. E-mail: leo@eco.meijigakuin.ac.jp}
\footnotetext[7]{Department of Mathematics, Faculty of Engineering, Osaka Institute of Technology. Omiya, Asahi-ku, Osaka 535-8585, Japan. E-mail: YHK03302@nifty.ne.jp}

This manuscript is the one which we submitted to Crelle Journal in September 2001. The first author talked on this subject at Oberwolfach "Theory of Riemann Zeta and Allied Functions" at 20.09.2001.

\section{Our Result}
This paper is a sequel to our previous paper \cite{CM}. Let $a(\geq 2)$ be a fixed natural number. For every residue class $l\pmod 4$, $l=0,1,2,3$, we consider the asymptotic behavior of the cardinality of the set 
$$Q_a(x;4,l):=\bigl\{p\leq x\ ;\ \sharp\langle a\pmod p\rangle\equiv l\pmod 4 \bigr\},$$
where $\sharp\langle a\pmod p\rangle$ denotes the order of the class $a\pmod p$ in $({\bf Z}/p{\bf Z})^\times$, the set of all invertible residue classes modulo a prime $p$. 

In our previous paper, we proved
\begin{thm}\label{th:prev}
We assume $a$ is a square free positive integer with $a\geq 3$. 

\noindent{\rm (I)} We have, for $l=0,2$,  
$$\sharp Q_a(x;4,l)=\frac{1}{3}\li x+O\left(\frac{x}{\log x\log\log x}\right).$$

\noindent{\rm (II)} We assume the Generalized Riemann Hypothesis (GRH) and further assume $a\equiv 1\pmod 4$. Then, for $l=1,3$, we have 
$$\sharp Q_a(x;4,l)=\frac{1}{6}\li x+O\left(\frac{x}{\log x\log\log x}\right).$$
\end{thm}
Although $\sharp Q_a(x;4,0)$ and $\sharp Q_a(x;4,2)$ are not difficult to study, yet the distributions of $Q_a(x;4,1)$ and $Q_a(x;4,3)$ are rather complicated. In the above theorem, we treated only the simplest case, but numerical examples show that, when $l=1,3$ , natural densities of $\sharp Q_a(x;4,l)$, with $a$ varies, are not always $1/6$.

In this paper, we remove those conditions on $a$ and will prove the following much more general result:
\begin{thm}\label{th:main}
We assume $a\in{\bf N}$ is not a perfect $h$-th power with $h\geq2$, and put 
$$a=a_1a_2^2,\quad a_1:\mbox{ square free.}$$
When $a_1\equiv 2\pmod 4$, we define $a_1'$ by 
$$a_1=2a_1'.$$
We assume GRH. And we define an absolute constant $C$ by 
\begin{equation}\label{eq:defC}
C:=\prod_{{\scr p\equiv3\,(\!\bmod4)}\atop{\scr p:\mbox{\scrs\rm prime}}}
\left(1-\frac{2p}{(p^2+1)(p-1)}\right).
\end{equation}
Then, for $l=1,3$, we have an asymptotic formula 
$$\sharp Q_a(x;4,l)=\delta_l\,\li x+O\left(\frac{x}{\log x\log\log x}\right),$$
and the leading coefficients $\delta_l$ ($l=1,3$) are given by the following way:

\noindent{\rm(I)} If $a_1\equiv 1,3\pmod 4$, then $\delta_1=\delta_3=1/6$.

\noindent{\rm(II)} When $a_1\equiv 2\pmod 4$, 

{\rm(i)} If $a_1'=1$, i.e. $a=2\cdot(\mbox{\rm a square number})$, then
$$\delta_1=\frac{7}{48}-\frac{C}{8},\quad \delta_3=\frac{7}{48}+\frac{C}{8}.$$

{\rm(ii)} If $a_1'\equiv1\pmod4$ with $a'_1>1$, then 

\ {\rm(ii-1)} if $a_1'$ has a prime divisor $p$ with $p\equiv1\pmod 4$, then $\delta_1=\delta_3=1/6$,

\ {\rm(ii-2)} if all prime divisors $p$ of $a_1'$ satisfy $p\equiv 3\pmod 4$, then 
\begin{eqnarray*}
\delta_1&=&\frac{1}{6}-\frac{C}{8}\prod_{p|a_1'}\left(\frac{-2p}{p^3-p^2-p-1}\right),\\
\delta_3&=&\frac{1}{6}+\frac{C}{8}\prod_{p|a_1'}\left(\frac{-2p}{p^3-p^2-p-1}\right).
\end{eqnarray*}

{\rm(iii)} If $a_1'\equiv 3\pmod 4$, then 

\ {\rm(iii-1)} if $a_1'$ has a prime divisor $p$ with $p\equiv1\pmod4$, then $\delta_1=\delta_3=1/6$,

\ {\rm(iii-2)} if all prime divisors $p$ of $a_1'$ satisfy $p\equiv 3\pmod 4$, then 
\begin{eqnarray*}
\delta_1&=&\frac{1}{6}+\frac{C}{8}\prod_{p|a_1'}\left(\frac{-2p}{p^3-p^2-p-1}\right),\\
\delta_3&=&\frac{1}{6}-\frac{C}{8}\prod_{p|a_1'}\left(\frac{-2p}{p^3-p^2-p-1}\right).
\end{eqnarray*}
\end{thm}

For our results, see also [ 3 ] and [ 4 ].

It seems an interesting phenomenon that, in (II)--(ii) and --(iii), the
densities $\delta_1$ and $\delta_3$ are controled by whether $a_1'$ has a prime factor $p$ with $p\equiv1 \pmod 4$ or not. Moreover, we can check easily that, in all cases, we have a {\it mysterious} inequality
$$\delta_1 \leq \delta_3.$$

In this paper, we limited our arguments to $\sharp Q_a(x;4,l)$ for $l=1,3$,  because these are more interesting than the other cases. Actually, we can prove {\it unconditionally} that
$$\begin{array}{ll}
\delta_0 = \delta_2 = 1/3, & \quad\mbox{if } a_1\ne2,\\
\delta_0=5/12 \mbox{ and }\delta_2=7/24, & \quad\mbox{if }a_1=2.
\end{array}$$

Throughout this paper, $p$ denotes an odd prime number, and for a natural number $n$, 
$$\zeta_n=\exp\left(\frac{2\pi i}{n}\right).$$
For an integer $l\geq0$, and for a natural number $f\geq1$, we define 
\begin{eqnarray*}
k&=&k(l,f)=2^f(4l+1),\\
k'&=&k'(l,f)=2^f(4l+3),
\end{eqnarray*}
and 
\begin{eqnarray*}
k_0&=&\prod_{p|k}p\quad(\mbox{the core of }k),\\
k'_0&=&\prod_{p|k'}p\quad(\mbox{the core of }k').
\end{eqnarray*}
For a square free integer $n\geq1$, and $d|k_0$, we construct extension fields
\begin{eqnarray*}
G_{k,n,d}&=&{\bf Q}(a^{1/kn}, \zeta_n, \zeta_{kd}),\\
\tilde G_{k,n,d}&=&G_{k,n,d}(\zeta_{2^{f+2}}).
\end{eqnarray*}
Furthermore, for $l=1,3$, let $\sigma_l$ be the automorphisms of ${\bf Q}(\zeta_{2^{f+2}})$ over ${\bf Q}$, which are defined by
\begin{eqnarray*}
\sigma_1:\zeta_{2^{f+2}}&\mapsto& \zeta_{2^{f+2}}^{\ 1+2^f},\\
\sigma_3:\zeta_{2^{f+2}}&\mapsto& \zeta_{2^{f+2}}^{\ 1+3\cdot 2^f}.
\end{eqnarray*}
Now we define the number $c^{(l)}(k,n,d)$ ($l=1,3$) by
$$c^{(l)}(k,n,d)=\left\{\begin{array}{ll}
1, & \mbox{if the map $\sigma_l$ can be extended\footnotemark[1] into the }\\
   & \mbox{automorphism $\sigma_l^\ast\in\Aut(\tilde G_{k,n,d}/G_{k,n,d})$,}\\[0.2cm]
0, & \mbox{otherwise.}
\end{array}\right.$$
\footnotetext[1]{``$\sigma_l^\ast\in\Aut(\tilde G_{k,n,d}/G_{k,n,d})$ is the extension of $\sigma_l$'' means that ``$\sigma_l^\ast|_{{\bf Q}(\zeta_{2^{f+2}})}=\sigma_l$.''}
Under the above notations, in \cite{CM} we obtained the following results: for $l=1,3$, under GRH, we have an asymptotic formula
$$\sharp Q_a(x;4,l)=\delta_l\,\li x+O\left(\frac{x}{\log x\log\log x}\right),$$
and the coefficients $\delta_1$ and $\delta_3$ are given by 
\begin{eqnarray}
\delta_1&=&\sum_{f\geq1}\sum_{l\geq0}
   \frac{k_0}{\varphi(k_0)}\sum_{d|k_0}\frac{\mu(d)}{d}\sum_{n}
   \frac{\mu(n)c^{(1)}(k,n,d)}{[\tilde G_{k,n,d}:{\bf Q}]}\nonumber\\
&+&\sum_{f\geq1}\sum_{l\geq0}
   \frac{k'_0}{\varphi(k'_0)}\sum_{d|k'_0}\frac{\mu(d)}{d}\sum_{n}
   \frac{\mu(n)c^{(3)}(k',n,d)}{[\tilde G_{k',n,d}:{\bf Q}]},
\label{eq:1-2}\\
\delta_3&=&\sum_{f\geq1}\sum_{l\geq0}
   \frac{k_0}{\varphi(k_0)}\sum_{d|k_0}\frac{\mu(d)}{d}\sum_{n}
   \frac{\mu(n)c^{(3)}(k,n,d)}{[\tilde G_{k,n,d}:{\bf Q}]}\nonumber\\
&+&\sum_{f\geq1}\sum_{l\geq0}
   \frac{k'_0}{\varphi(k'_0)}\sum_{d|k'_0}\frac{\mu(d)}{d}\sum_{n}
   \frac{\mu(n)c^{(1)}(k',n,d)}{[\tilde G_{k',n,d}:{\bf Q}]}.
\label{eq:1-3}
\end{eqnarray}
In our proof of Theorem \ref{th:prev} (see \cite[Theorem 1.2]{CM}), we compared the coefficients (\ref{eq:1-2}) and (\ref{eq:1-3}). And for a square free $a$ with $a\equiv 1\pmod4$, we can prove that, for any $k,k',n,d$, 
\begin{eqnarray*}
c^{(1)}(k,n,d)&=&c^{(3)}(k,n,d)\\
c^{(1)}(k',n,d)&=&c^{(3)}(k',n,d).
\end{eqnarray*}
Then the above expressions (\ref{eq:1-2}) and (\ref{eq:1-3}) give the same number --- this is the key idea of our previous work. This method is available only to prove $\delta_1=\delta_3$, but numerical examples show that the equality $\delta_1=\delta_3$ does not always hold ({\it cf.} Tables 2.1 and 2.2 in \S 2, and for example $a=2$). In this paper, we calculate the infinite sums (\ref{eq:1-2}) and (\ref{eq:1-3}) directly. 

We calculate the extension degree $[\tilde G_{k,n,d}:{\bf Q}]$ in \S 3, determine the value of the coefficients $c^{(l)}(k,n,d)$, $l=1,3$ in \S 4. And we will prove Theorem \ref{th:main} in \S 5. Preceding our proof, in order to clarify our results, we mention some numerical examples in \S 2. 

In what follows, $\lcm{x,y}$ denotes the least common multiple of $x$ and $y$, 
$$\underline{x}=\mbox{the odd part of }x,$$
and $\psum_n$ means the sum over ``square free'' numbers $n$'s. 
%
%
\section{Numerical Examples}
We will compare the theoretical density $\delta_l$ with the experimental density $\sharp Q_a(x;4,l)/\pi(x)$ for various $a$'s. Here we take $x=10^7$. And according to our results, we omit the cases $l=0,2$ ({\it cf.} also Tables 5.1 -- 5.4 in \cite{CM}). As for the number $C$ defined in (\ref{eq:defC}), we use a rough approximate value 
$$C\approx 0.64365.$$

\medskip
\newcommand{\lw}[1]{\smash{\lower1.7ex\hbox{#1}}}
{\footnotesize 
\begin{center}
\renewcommand{\arraystretch}{1.2}
\begin{tabular}{|c|c|l||c|c|}
\hline
\multicolumn{2}{|c|}{type of $a_1$} & \multicolumn{1}{|c||}{$a$} & the theoretical density $\delta_1$ & $\sharp Q_a(10^7;4,l)/\pi(10^7)$\\
\hline\hline
\multicolumn{2}{|c|}{     } & 5 & 1/6 & 0.166771\\
\multicolumn{2}{|c|}{$a_1\equiv 1\pmod4$} & $33=3\cdot11$ & 1/6 & 0.166991\\
\multicolumn{2}{|c|}{     } & $45=5\cdot 3^2$ & 1/6 & 0.167141\\
\hline
 & \lw{$a_1'=1$} & 2 & $7/48-C/8\approx 0.06538$ & 0.065425\\
 & & $50=2\cdot5^2$ & $7/48-C/8$ & 0.065351\\
\cline{2-5}
 &  & $10=2\cdot5$ & 1/6 & 0.166644\\
\lw{$a_1\equiv 2\pmod4$} & $a_1'\equiv1\pmod4$ & $42=2\cdot3\cdot7$ & $1/6-3C/1144\approx 0.16498$ & 0.165129\\
 & & $210=2\cdot3\cdot5\cdot7$ & $1/6$ & 0.166878\\
\cline{2-5}
 & \lw{} & $6=2\cdot3$ & $1/6-3C/56\approx0.13219$ & 0.132179\\
 & $a_1'\equiv3\pmod4$ & $14=2\cdot7$ & $1/6-7C/1144\approx0.16273$ & 0.162875\\
 & & $30=2\cdot3\cdot5$ & 1/6 & 0.166354\\
\hline
\multicolumn{2}{|c|}{     } & 11 & 1/6 & 0.166531\\
\multicolumn{2}{|c|}{{$a_1\equiv3\pmod4$}} & $55=5\cdot11$ & 1/6 & 0.166875\\
\multicolumn{2}{|c|}{     } & $75=3\cdot5^2$ & 1/6 & 0.166372\\
\hline
\end{tabular}
\renewcommand{\arraystretch}{1.0}

\smallskip
{\bf Table 2.1.}\label{tb:2-1} The densities of $Q_a(x;4,1)$. 
\end{center}
}

\medskip
{\footnotesize 
\begin{center}
\renewcommand{\arraystretch}{1.2}
\begin{tabular}{|c|c|l||c|c|}
\hline
\multicolumn{2}{|c|}{type of $a_1$} & \multicolumn{1}{|c||}{$a$} & the theoretical density $\delta_1$ & $\sharp Q_a(10^7;4,l)/\pi(10^7)$\\
\hline\hline
\multicolumn{2}{|c|}{\lw{}} & 5 & 1/6 & 0.166810\\
\multicolumn{2}{|c|}{$a_1\equiv 1\pmod4$} & $33=3\cdot11$ & 1/6 & 0.166274\\
\multicolumn{2}{|c|}{     } & $45=5\cdot 3^2$ & 1/6 & 0.166324\\
\hline
 & \lw{$a_1'=1$} & 2 & $7/48+C/8\approx 0.22629$ & 0.226407\\
 & & $50=2\cdot5^2$ & $7/48+C/8$ & 0.226345\\
\cline{2-5}
 & \lw{} & $10=2\cdot5$ & 1/6 & 0.166522\\
\lw{$a_1\equiv 2\pmod4$} & $a_1'\equiv1\pmod4$ & $42=2\cdot3\cdot7$ & $1/6+3C/1144\approx 0.16835$ & 0.168277\\
 & & $210=2\cdot3\cdot5\cdot7$ & 1/6 & 0.166490\\
\cline{2-5}
 & \lw{} & $6=2\cdot3$ & $1/6+3C/56\approx0.20115$ & 0.201471\\
 & $a_1'\equiv3\pmod4$ & $14=2\cdot7$ & $1/6+7C/1144\approx0.17061$ & 0.170289\\
 & & $30=2\cdot3\cdot5$ & 1/6 & 0.166991\\
\hline
\multicolumn{2}{|c|}{     } & 11 & 1/6 & 0.166766\\
\multicolumn{2}{|c|}{{$a_1\equiv3\pmod4$}} & $55=5\cdot11$ & 1/6 & 0.166691\\
\multicolumn{2}{|c|}{     } & $75=3\cdot5^2$ & 1/6 & 0.166896\\
\hline
\end{tabular}
\renewcommand{\arraystretch}{1.0}

\smallskip
{\bf Table 2.2.}\label{tb:2-2} The densities of $Q_a(x;4,3)$. 
\end{center}
}
%
%
\section{The Extension Degree $[\tilde G_{k,n,d}:{\bf Q}]$}
In order to calculate the extension degree $[\tilde G_{k,n,d}:{\bf Q}]$ which appears in (\ref{eq:1-2}) and (\ref{eq:1-3}), we need a few lemmas. 
\begin{lem}\label{lem:3-1}
Let $u$ be a natural number. The maximal normal subfield which is contained in ${\bf Q}(a^{1/u})$ is
$$\left\{\begin{array}{ll}
{\bf Q}, & \mbox{if $u$ is odd,}\\
{\bf Q}(\sqrt{a_1}), & \mbox{if $u$ is even.}
\end{array}
\right.$$
\end{lem}
We omit the proof.
\begin{lem}\label{lem:3-2}
The minimal cyclotomic field which contains ${\bf Q}(\sqrt{a_1})$ is 
$$\left\{\begin{array}{ll}
{\bf Q}(\zeta_{a_1}), & \mbox{if }a_1\equiv1\pmod4\\
{\bf Q}(\zeta_{4a_1}), & \mbox{if }a_1\equiv2,3\pmod4
\end{array}
\right.$$
\end{lem}
\prf Moree \cite[Lemma 1]{Mo}. \qed
\begin{lem}\label{lem:3-3}
Let $v_1$, $v_2$ be two natural numbers. When ${\bf Q}(\zeta_{v_1})\supset{\bf Q}(\zeta_{v_2})$, then $v_2|v_1$ or $v_2|2v_1$ with $v_2$: even and $v_1$: odd.
\end{lem}
\prf We put
$$V=(v_1,v_2)\quad \mbox{and} \quad\begin{array}{ll}
v_1=Vv_1'\\
v_2=Vv_2',
\end{array}$$
then 
$${\bf Q}(\zeta_{v_1})\cap{\bf Q}(\zeta_{v_2})={\bf Q}(\zeta_{V})={\bf Q}(\zeta_{Vv_2'}).$$
Thus 
\begin{equation}\label{eq:3-1}
\varphi(V)=\varphi(Vv_2').
\end{equation}
When $v_2$ is odd, $v_2'$ must be one, i.e. $v_2|v_1$. If $v_2$ is odd and $v_2\ne1$, then (\ref{eq:3-1}) yields $V$ is odd and $v_2'=2$. This shows that $v_2=2V$ and $v_1=Vv_1'$ with $v_1'$: odd. \qed
\begin{prop}\label{prop:3-1}
Let $u$ and $v$ be natural numbers, and 
$$D_{u,v}=[{\bf Q}(a^{1/u},\zeta_{v}):{\bf Q}].$$
Then we have

\noindent{\rm(1)} When $u$ is odd, $D_{u,v}=u\varphi(v)$. 

\noindent{\rm(2)} When $u$ is even, $D_{u,v}=u\varphi(v)$ or $u\varphi(v)/2$ and the latter happens, if and only if, one of {\rm(i), (ii), (iii)} is satisfied:
$$\begin{array}{llll}
{\rm(i)} & a_1\equiv1\pmod4 & \mbox{\rm and} & a_1|v\\
{\rm(ii)} & a_1\equiv2\pmod4 & \mbox{\rm and} & 4a_1|v\\
{\rm(iii)} & a_1\equiv3\pmod4 & \mbox{\rm and} & 4a_1|v.
\end{array}$$
\end{prop}
\prf From the assumption of $a$, we have 
$$[{\bf Q}(a^{1/u}):{\bf Q}]=u,$$
then 
\begin{eqnarray*}
D_{u,v}&=&\frac{[{\bf Q}(a^{1/u}):{\bf Q}][{\bf Q}(\zeta_{v}):{\bf Q}]}
           {[{\bf Q}(a^{1/u})\cap{\bf Q}(\zeta_{v}):{\bf Q}]}\\
&=&\frac{1}{[{\bf Q}(a^{1/u})\cap{\bf Q}(\zeta_{v}):{\bf Q}]}\cdot u\varphi(v).
\end{eqnarray*}
The field ${\bf Q}(a^{1/u})\cap{\bf Q}(\zeta_{v})$ is a normal extension of ${\bf Q}$, then, by Lemma \ref{lem:3-1}, 
$${\bf Q}(a^{1/u})\cap{\bf Q}(\zeta_{v})={\bf Q} \quad\mbox{or}\quad {\bf Q}(\sqrt{a_1}).$$
The latter happens, if and only if, ${\bf Q}(\sqrt{a})\subset{\bf Q}(\zeta_{v})$. Thus, when $a_1\equiv1\pmod4$, Lemma \ref{lem:3-2} implies that ${\bf Q}(\zeta_{v})\supset {\bf Q}(\zeta_{a_1})$, and Lemma \ref{lem:3-3} shows $a_1|v$. 

When $a_1\equiv2,3\pmod4$, Lemma \ref{lem:3-2} implies ${\bf Q}(\zeta_{v})\supset {\bf Q}(\zeta_{4a_1})$, and Lemma \ref{lem:3-3} shows again $4a_1|v$. \qed
\begin{cor}[{The Values of $[\tilde G_{k,n,d}:{\bf Q}]$}]
\label{cor:3-1}
When one of {\rm (i), (ii), (iii)} is satisfied, then
$$[\tilde G_{k,n,d}:{\bf Q}]=\frac{1}{2}nk\cdot 2^{f+1}\varphi(\lcm{\ul{n},\ul{k}\ul{d}}),$$
$$\begin{array}{llll}
{\rm(i)} & a_1\equiv1\pmod4 & \mbox{\rm and} & a_1|\lcm{\ul{n},\ul{k}\ul{d}},\\
{\rm(ii)} & a_1\equiv2\pmod4 & \mbox{\rm and} & 
     a_1'|\lcm{\ul{n},\ul{k}\ul{d}},\ \mbox{\rm where } a_1'=a_1/2,\\
{\rm(iii)} & a_1\equiv3\pmod4 & \mbox{\rm and} & a_1|\lcm{\ul{n},\ul{k}\ul{d}}.
\end{array}$$
Otherwise,
$$[\tilde G_{k,n,d}:{\bf Q}]=nk\cdot 2^{f+1}\varphi(\lcm{\ul{n},\ul{k}\ul{d}}).$$
When we exchange $k$ for $k'$, then we get the values of $[\tilde G_{k',n,d}:{\bf Q}]$.
\end{cor}
\prf From the definition of $\tilde G_{k,n,d}$, 
$$\tilde G_{k,n,d}={\bf Q}(a^{1/nk},\zeta_{n},\zeta_{kd},\zeta_{2^{f+2}})$$
and
$$\lcm{n,kd,2^{f+2}}=2^{f+2}\lcm{\ul{n},\ul{k}\ul{d}}.\ \qed$$
\begin{cor}\label{cor:3-2}
We assume $d$ is odd, then the extension degree $[\tilde G_{k,n,d}:G_{k,n,d}]=4$ or $2$, and the latter case happens, if and only if, one of {\rm (i)$'$, (ii)$'$, (iii)$'$} is satisfied:
$$\begin{array}{lllll}
{\rm(i)'} & a_1\equiv2\pmod4, & d:\mbox{\rm odd}, & a_1'|\lcm{\ul{n},\ul{k}\ul{d}} & \mbox{\rm and } f=1,\\
{\rm(ii)'} & a_1\equiv2\pmod4, & d:\mbox{\rm odd}, & a_1'|\lcm{\ul{n},\ul{k}\ul{d}} & \mbox{\rm and } f=2,\\
{\rm(iii)'} & a_1\equiv3\pmod4, & d:\mbox{\rm odd},  & a_1|\lcm{\ul{n},\ul{k}\ul{d}} & \mbox{\rm and } f=1.
\end{array}$$
\end{cor}
\prf $[\tilde G_{k,n,d}:G_{k,n,d}]=2$ happens, if and only if, 
$$[\tilde G_{k,n,d}:{\bf Q}]=\frac12 nk\cdot 2^{f+1}\varphi(\lcm{\ul{n},\ul{k}\ul{d}})$$
and 
\begin{equation}\label{eq:3-2}
[G_{k,n,d}:{\bf Q}]=nk\cdot 2^{f-1}\varphi(\lcm{\ul{n},\ul{k}\ul{d}}).
\end{equation}
Proposition \ref{prop:3-1} shows that (\ref{eq:3-2}) happens, if and only if, one of the following {\rm (i)$''$, (ii)$''$, (iii)$''$} is satisfied:
$$\begin{array}{lll}
{\rm(i)''} & a_1\equiv1\pmod4, & a_1\nmid\lcm{n,kd},\\
{\rm(ii)''} & a_1\equiv2\pmod4, & 4a_1'\nmid\lcm{n,kd},\\
{\rm(iii)''} & a_1\equiv3\pmod4, & 4a_1\nmid\lcm{n,kd}.
\end{array}$$
Combining {\rm (i)$''$, (ii)$''$, (iii)$''$} with {\rm (i), (ii), (iii)} of Corollary \ref{cor:3-1}, we can verify {\rm (i)$'$, (ii)$'$, (iii)$'$} easily. \qed
%
%
\section{The Coefficients $c^{(l)}(k,n,d)$}
In this section, we consider the value of the coefficient $c^{(l)}(k,n,d)$, which appeared in (\ref{eq:1-2}) and (\ref{eq:1-3}). 
\begin{lem}\label{lem:4-1}
If $d$ is even, then for any $k,n$, 
$$c^{(1)}(k,n,d)=c^{(3)}(k,n,d)=0,$$
and the same for $k'$. 
\end{lem}
\prf Let $\sigma_1^\ast$ be the extension of $\sigma_1$. Then
\begin{eqnarray*}
\sigma_1^\ast(\zeta_{2^{f+1}})&=&\sigma_1^\ast(\zeta_{2^{f+2}}^2)
  =\{(\zeta_{2^{f+2}})^{1+2^f}\}^2\\
 &=&\zeta_{2^{f+2}}^{2+2^{f+1}}=-\zeta_{2^{f+1}},
\end{eqnarray*}
and this contradicts the condition $\sigma_1^\ast|_{G_{k,n,d}}=\id$. This proves $c^{(1)}(k,n,d)=0$, and similarly  $c^{(3)}(k,n,d)=0$. \qed
\begin{lem}\label{lem:4-2}
If $d$ is odd and $[\tilde G_{k,n,d}:G_{k,n,d}]=4$, then $c^{(1)}(k,n,d)=c^{(3)}(k,n,d)=1$, and the same for $k'$. 
\end{lem}
\prf Both two extensions $\tilde G_{k,n,d}/G_{k,n,d}$ and ${\bf Q}(\zeta_{2^{f+2}})/{\bf Q}(\zeta_{2^{f}})$ are Galois extensions. And when $[\tilde G_{k,n,d}:G_{k,n,d}]=4$, 
$$\Gal(\tilde G_{k,n,d}/G_{k,n,d})\cong\Gal({\bf Q}(\zeta_{2^{f+2}})/{\bf Q}(\zeta_{2^{f}})).$$
So we can extend $\sigma_1, \sigma_3\in\Gal({\bf Q}(\zeta_{2^{f+2}})/{\bf Q}(\zeta_{2^{f}}))$ uniquely to $\sigma_1^\ast, \sigma_3^\ast\in\Gal(\tilde G_{k,n,d}/G_{k,n,d})$, and we have 
$$c^{(1)}(k,n,d)=c^{(3)}(k,n,d)=1. \ \qed
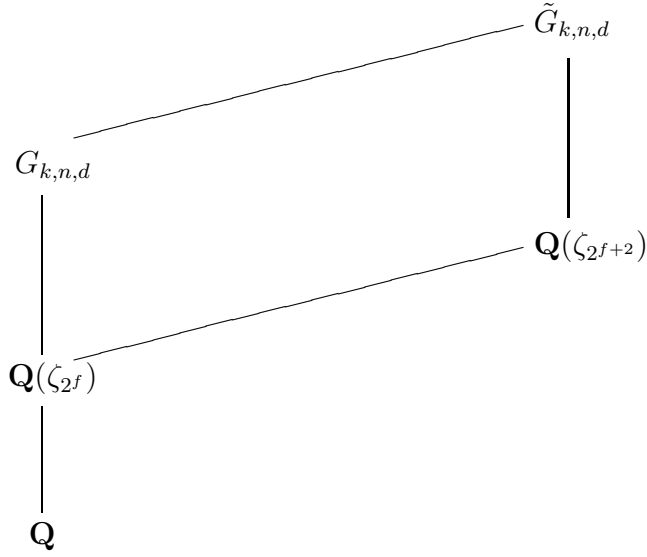
$$
\begin{center}
\begin{picture}(270,250)
\put(-4,20){\makebox(20,14)[c]{${\bf Q}$}}
\put(0,80){\makebox(20,14)[c]{${\bf Q}(\zeta_{2^f})$}}
\put(0,160){\makebox(20,14)[c]{$G_{k,n,d}$}}
%
\put(192,130){\makebox(20,14)[l]{${\bf Q}(\zeta_{2^{f+2}})$}}
\put(192,215){\makebox(20,14)[l]{$\tilde G_{k,n,d}$}}
%
\put(6,36){\line(0,1){40}}
\put(6,96){\line(0,1){60}}
%
\put(205,148){\line(0,1){60}}
%
\put(18,94){\line(4,1){170}}
\put(18,178){\line(4,1){170}}
%
\end{picture}
\end{center}
\begin{lem}\label{lem:4-3}
If $d$ is odd and $f\geq2$, then, for any $l,n$,
$$c^{(1)}(k,n,d)=c^{(3)}(k,n,d),$$
and the same for $k'$. 
\end{lem}
\prf When the extension $\sigma_1^\ast$ exists, then, since $\sigma_1^\ast\circ\sigma_1^\ast\circ\sigma_1^\ast=\sigma_3^\ast$, $\sigma_3^\ast$ exists, and vice versa. \qed

\medskip
Now we consider the case ``$d$: odd and $[\tilde G_{k,n,d}:G_{k,n,d}]=2$.''

\medskip
\noindent \ul{\bf Case (i)$'$ of Corollary \ref{cor:3-2}}

\medskip
\nopagebreak
\noindent $f=1$, $d$: odd, $a_1\equiv2\pmod4$ and $a_1'|\lcm{\ul{n},\ul{k}d}$.
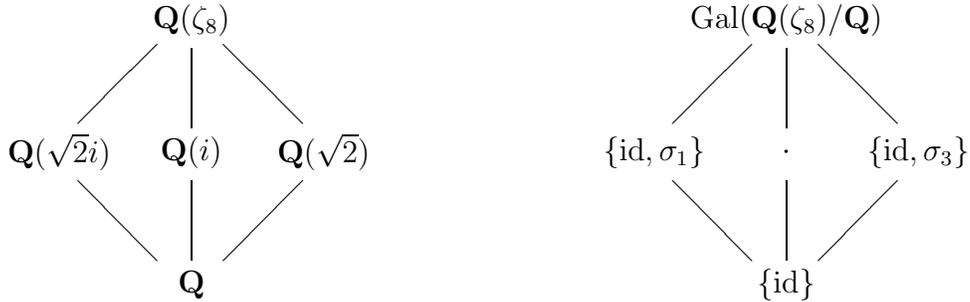
\begin{center}
\begin{picture}(340,120)
\put(0,50){\makebox(20,14)[c]{${\bf Q}(\sqrt2i)$}}
\put(50,0){\makebox(20,14)[c]{${\bf Q}$}}
\put(50,50){\makebox(20,14)[c]{${\bf Q}(i)$}}
\put(50,100){\makebox(20,14)[c]{${\bf Q}(\zeta_{8})$}}
\put(100,50){\makebox(20,14)[c]{${\bf Q}(\sqrt2)$}}
\put(60,16){\line(0,1){30}}
\put(60,66){\line(0,1){30}}
\put(17,68){\line(1,1){30}}
\put(72,16){\line(1,1){30}}
\put(17,46){\line(1,-1){30}}
\put(72,98){\line(1,-1){30}}
\put(200,50){\makebox(70,14)[c]{$\{\id, \sigma_1\}$}}
\put(250,0){\makebox(70,14)[c]{$\{\id\}$}}
\put(250,50){\makebox(70,14)[c]{$\cdot$}}
\put(250,100){\makebox(70,14)[c]{$\Gal({\bf Q}(\zeta_8)/{\bf Q})$}}
\put(300,50){\makebox(70,14)[c]{$\{\id,\sigma_3\}$}}
\put(285,16){\line(0,1){30}}
\put(285,66){\line(0,1){30}}
\put(242,68){\line(1,1){30}}
\put(297,16){\line(1,1){30}}
\put(242,46){\line(1,-1){30}}
\put(297,98){\line(1,-1){30}}
\end{picture}
\end{center}
We define an automorphism $\tau$ by $\Gal(\tilde G_{k,n,d}/G_{k,n,d})=\{\id_{\tilde G_{k,n,d}},\tau\}$. We remark that
\begin{eqnarray}
\tau(\sqrt2i)=\sqrt2i &\Leftrightarrow& c^{(1)}(k,n,d)=1,\ c^{(3)}(k,n,d)=0
\label{eq:4-1}\\
\tau(i)=i &\Leftrightarrow& c^{(1)}(k,n,d)=c^{(3)}(k,n,d)=0
\label{eq:4-2}\\
\tau(\sqrt2)=\sqrt2 &\Leftrightarrow& c^{(1)}(k,n,d)=0,\ c^{(3)}(k,n,d)=1.
\label{eq:4-3}
\end{eqnarray}
Here we further assume that \ul{$a_1'\equiv1\pmod4$}. Then, by Lemma \ref{lem:3-2}, ${\bf Q}(\sqrt{a_1'})\subset{\bf Q}(\zeta_{\lcm{\ul{n},kd}})$ and 
$$\tau(\sqrt{a_1'})=\sqrt{a_1'}.$$
Since $\sqrt{a_1}\in G_{k,n,d}$ a priori, $\tau(\sqrt{a_1})=\sqrt{a_1}$. Thus 
$$\tau(\sqrt{2})=\sqrt{2}$$
and by (\ref{eq:4-3}), we have 
\begin{equation}
c^{(1)}(k,n,d)=0\quad \mbox{and}\quad c^{(3)}(k,n,d)=1.
\label{eq:4-4}\\
\end{equation}
When \ul{$a_1'\equiv3\pmod4$}, we consider the field
$$G'_{k,n,d}=G_{k,n,d}(\zeta_4).$$
It is easily seen that $G'_{k,n,d}\ne G_{k,n,d}$, this shows that
\begin{equation}
i\not\in G_{k,n,d}.
\label{eq:4-5}
\end{equation}
Now we will prove $\sqrt{2}i\in G_{k,n,d}$. The condition ``$a_1'\equiv3\pmod4$ and $4\nmid\lcm{n,kd}$'' implies 
$$\sqrt{a_1'}\not\in{\bf Q}(\zeta_{\lcm{n,kd}}),$$
and ``$a_1'\equiv3\pmod4$ and $4|\lcm{n,kd,4}$'' implies 
$$\sqrt{a_1'}\in{\bf Q}(\zeta_{\lcm{n,kd}},\zeta_4).$$
Since $[G'_{k,n,d}:{\bf Q}]=2[G_{k,n,d}:{\bf Q}]$, we have 
$$G_{k,n,d}\not\supset{\bf Q}(\zeta_{\lcm{n,kd}},\zeta_4)={\bf Q}(\zeta_{\lcm{n,kd}},\sqrt{a_1'}),$$
therefore
$$\sqrt{a_1'}\not\in G_{k,n,d}.$$
We have $\sqrt{2a_1'}\in G_{k,n,d}$ a priori, then
$$\sqrt2\not\in G_{k,n,d}.$$
Combining with (\ref{eq:4-5}), we can conclude $\sqrt{2}i\in G_{k,n,d}$, and 
\begin{equation}
c^{(1)}(k,n,d)=1\quad \mbox{and}\quad c^{(3)}(k,n,d)=0.
\label{eq:4-6}
\end{equation}

\medskip
\noindent \ul{\bf Case (ii)$'$ of Corollary \ref{cor:3-2}}

\medskip
\nopagebreak
\noindent $f=2$, $d$: odd, $a_1\equiv2\pmod4$ and $a_1'|\lcm{\ul{n},\ul{k}d}$.

\medskip
\noindent In this case, from Lemma \ref{lem:4-3}, we have already 
$$c^{(1)}(k,n,d)=c^{(3)}(k,n,d).$$
Since $[\tilde G_{k,n,d}:G_{k,n,d}]=2$, $\sharp \Gal(\tilde G_{k,n,d}/G_{k,n,d})=2$, then 
``$c^{(1)}(k,n,d)=c^{(3)}(k,n,d)=1$'' is impossible. Therefore
\begin{equation}
c^{(1)}(k,n,d)=c^{(3)}(k,n,d)=0.
\label{eq:4-7}
\end{equation}

\medskip
\noindent \ul{\bf Case (iii)$'$ of Corollary \ref{cor:3-2}}

\medskip
\nopagebreak
\noindent $f=1$, $d$: odd, $a_1\equiv3\pmod4$ and $a_1|\lcm{\ul{n},\ul{k}d}$.

\medskip
\noindent We consider again the field $G'_{k,n,d}=G_{k,n,d}(\zeta_4)$. A simple application of Proposition \ref{prop:3-1} (2) shows that 
$$[G'_{k,n,d}:{\bf Q}]=\frac12 nk\varphi(4\lcm{\ul{n},\ul{k}d}),$$
this means $G'_{k,n,d}=G_{k,n,d}$, i.e. 
$$\zeta_4\in G_{k,n,d}.$$
Therefore 
\begin{equation}
c^{(1)}(k,n,d)=c^{(3)}(k,n,d)=0.
\label{eq:4-8}
\end{equation}
It is clear that these arguments are true for $k'$ instead of $k$. 

Summing up the above results, we get the following table:
\begin{prop}\label{prop:4-1}
The value of the coefficients $c^{(l)}(k,n,d)$ and $c^{(l)}(k',n,d)$, $l=1,3$, are given as follows:

\noindent{\rm(1)} When $d$ is even, $c^{(l)}(k,n,d)=c^{(l)}(k',n,d)=0$.

\noindent{\rm(2)} When $d$ is odd, 
{\rm 
\begin{center}
\renewcommand{\arraystretch}{1.2}
\begin{tabular}{llccc}
 & & $c^{(1)}(k,n,d)$ & $c^{(3)}(k,n,d)$ & \\
\hline
(i) $a_1\equiv1\pmod4$ & & 1 & 1 & \\
\hline
(ii) $a_1\equiv2\pmod4$ & & & & \\
\quad (a) $a_1'\equiv1\pmod4$ & & & & \\
\qquad if $a_1'\nmid\lcm{\ul{n},\ul{k}d}$ & & 1 & 1 & \\
\qquad if $a_1'|\lcm{\ul{n},\ul{k}d}$, & $f=1$ & 0 & 1 & (\ref{eq:4-4})\\
                                       & $f=2$ & 0 & 0 & (\ref{eq:4-7})\\
                                    & $f\geq3$ & 1 & 1 & \\
\quad (b) $a_1'\equiv3\pmod4$ & & & & \\
\qquad if $a_1'\nmid\lcm{\ul{n},\ul{k}d}$ & & 1 & 1 & \\
\qquad if $a_1'|\lcm{\ul{n},\ul{k}d}$, & $f=1$ & 1 & 0 & (\ref{eq:4-6})\\
                                       & $f=2$ & 0 & 0 & (\ref{eq:4-7})\\
                                    & $f\geq3$ & 1 & 1 & \\
\hline
(iii) $a_1\equiv3\pmod4$ & & & & \\
\qquad if $a_1\nmid\lcm{\ul{n},\ul{k}d}$ & & 1 & 1 & \\
\qquad if $a_1|\lcm{\ul{n},\ul{k}d}$, & $f=1$ & 0 & 0 & (\ref{eq:4-8})\\
                                    & $f\geq2$ & 1 & 1 & \\
\end{tabular}
\renewcommand{\arraystretch}{1.0}
\end{center}
}
\noindent And the same results hold for $k'$.
\end{prop}
Here we remark that, except for the two cases (\ref{eq:4-4}) and (\ref{eq:4-6}), we have always $c^{(1)}(k,n,d)=c^{(3)}(k,n,d)$.
%
%
\section{Proof of Theorem \ref{th:main}}
In order to calculate the infinite sums (\ref{eq:1-2}) and (\ref{eq:1-3}), we first consider the following sums: let
\begin{eqnarray*}
I^{(1)}(f)&=&\sum_{l\geq0}\frac{k_0}{\varphi(k_0)}
  \sum_{\scriptstyle d|k_0 \atop \scriptstyle d:{\rm odd}}\frac{\mu(d)}{d}
  \psum_{n}\frac{\mu(n)}{nk\varphi(\lcm{\ul{n},\ul{k}d})},\\
I^{(3)}(f)&=&\sum_{l\geq0}\frac{k'_0}{\varphi(k'_0)}
  \sum_{\scriptstyle d|k'_0 \atop \scriptstyle d:{\rm odd}}\frac{\mu(d)}{d}
  \psum_{n}\frac{\mu(n)}{nk'\varphi(\lcm{\ul{n},\ul{k'}d})},
\end{eqnarray*}
and for square free integer $s\geq1$, let
\begin{eqnarray*}
J^{(1)}(f,s)&=&
  \sum_{\scriptstyle l\geq0 \atop\scriptstyle s|\lcm{\ul{n},\ul{k}d}}
  \frac{k_0}{\varphi(k_0)}
  \sum_{\scriptstyle d|k_0 \atop \scriptstyle d:{\rm odd}}\frac{\mu(d)}{d}
  \psum_{n}\frac{\mu(n)}{nk\varphi(\lcm{\ul{n},\ul{k}d})},\\
J^{(3)}(f,s)&=&
  \sum_{\scriptstyle l\geq0 \atop\scriptstyle s|\lcm{\ul{n},\ul{k'}d}}
  \frac{k'_0}{\varphi(k'_0)}
  \sum_{\scriptstyle d|k'_0 \atop \scriptstyle d:{\rm odd}}\frac{\mu(d)}{d}
  \psum_{n}\frac{\mu(n)}{nk'\varphi(\lcm{\ul{n},\ul{k'}d})}.
\end{eqnarray*}
We remark here that, when $s=1$, $J^{(l)}(f,1)=I^{(l)}(f)$ and the sums $J^{(1)}(f,s)$ and $J^{(3)}(f,s)$ are partial sums of $I^{(1)}(f)$ and $I^{(3)}(f)$, respectively. 

We can calculate these sums as follows:
\begin{lem}\label{lem:5-1}
For any $f\geq1$, 
$$I^{(1)}(f)+I^{(3)}(f)=\frac{1}{2^f}.$$
\end{lem}
\begin{lem}\label{lem:5-2}
For any $f\geq1$ and any square free $s>1$, 
$$J^{(1)}(f,s)+J^{(3)}(f,s)=0.$$
\end{lem}
\begin{lem}\label{lem:5-3}
We use the number $C$ defined by (\ref{eq:defC}). For any $f\geq1$, we have
{\rm
$$J^{(1)}(f,s)=\left\{\begin{array}{cl}
\Frac{1+C}{2^{f+1}}, & \mbox{if }s=1,\\[0.3cm]
0, & \mbox{if $s$ has a prime divisor $p$ with $p\equiv1\pmod4$},\\[0.2cm]
\Frac{C}{2^{f+1}}\prod_{p|s}\Frac{-2p}{p^3-p^2-p-1}, & \mbox{if all prime divisors $p$ of $s$ satisfy $p\equiv3\pmod4$,}
\end{array}\right.$$
$$J^{(3)}(f,s)=\left\{\begin{array}{cl}
\Frac{1-C}{2^{f+1}}, & \mbox{if }s=1,\\[0.3cm]
0, & \mbox{if $s$ has a prime divisor $p$ with $p\equiv1\pmod4$},\\[0.2cm]
\Frac{-C}{2^{f+1}}\prod_{p|s}\Frac{-2p}{p^3-p^2-p-1}, & \mbox{if all prime divisors $p$ of $s$ satisfy $p\equiv3\pmod4$.}
\end{array}\right.$$
}
\end{lem}

\medskip
\noindent{\bf Proof of Lemma \ref{lem:5-1}.}

\medskip
\noindent We put the primary decomposition of $k$ as
\begin{equation}
k=2^f(4l+1)=2^fq_1^{e_1}\cdots q_r^{e_r}.
\label{eq:5-1}
\end{equation}
Since $d$ is odd and $d|k_0$, we have
$$d=q_1^{\varepsilon_1}\cdots q_r^{\varepsilon_r}
\quad\mbox{with }\varepsilon_i=0 \mbox{ or }1.$$
For a square free $n$, 
$$\lcm{\ul{n},\ul{k}d}=
q_1^{e_1+\varepsilon_1}\cdots q_r^{e_r+\varepsilon_r}
\prod_{\scr p|\ul{n} \atop\scr p\ne q_i}p,$$
and
\begin{eqnarray}
\sum_{\scr d|k_0 \atop\scr d:\mbox{\scrs odd}}
\frac{\mu(d)}{d}\cdot\frac{1}{\varphi(\lcm{\ul{n},\ul{k}d})}&=&
\sum_{\varepsilon_i=0,1}
\frac{(-1)^{\varepsilon_1}\cdots(-1)^{\varepsilon_r}}
{q_1^{\varepsilon_1}\cdots q_r^{\varepsilon_r}}
\prod_{\scr p|\ul{n} \atop \scr p\ne q_i}\frac{1}{p-1}
\prod_{i=1}^{r}\frac{1}{(q_i-1)q_i^{e_i+\varepsilon_i-1}}\nonumber\\
&=&\prod_{\scr p|\ul{n} \atop \scr p\ne q_i}\frac{1}{p-1}
\prod_{i=1}^{r}\left(\frac{1}{(q_i-1)q_i^{e_i-1}}+\frac{-1}{(q_i-1)q_i^{e_i+1}}
\right)\nonumber\\
&=&\prod_{\scr p|\ul{n} \atop \scr p\ne q_i}\frac{1}{p-1}
\prod_{i=1}^{r}\frac{q_i+1}{q_i^{e_i+1}}.\label{eq:5-2}
\end{eqnarray}
Now
\begin{eqnarray}
\psum_{n}\frac{\mu(n)}{n}\prod_{\scr p|\ul{n} \atop \scr p\ne q_i}\frac{1}{p-1}
&=&\left(\psum_{n:\mbox{\scrs odd}}+\psum_{n:\mbox{\scrs even}}\right)
\frac{\mu(n)}{n}\prod_{\scr p|\ul{n} \atop \scr p\ne q_i}\frac{1}{p-1}
\nonumber\\
&=&\frac12\prod_{i=1}^{r}\left(1+\frac{-1}{q_i}\right)
\prod_{\scr p\ne q_i \atop \scr p:\mbox{\scrs odd}}\left(1+\frac{-1}{p(p-1)}\right).
\label{eq:5-3}
\end{eqnarray}
Consequently we have
\begin{eqnarray}
I^{(1)}(f)&=&\sum_{l\geq0}\frac{k_0}{k\varphi(k_0)}
\prod_{i=1}^{r}\frac{q_i+1}{q_i^{e_i+1}}\cdot\frac12
\prod_{i=1}^{r}\left(\frac{q_i-1}{q_i}\right)
\prod_{\scr p\ne q_i \atop \scr p:\mbox{\scrs odd}}
\left(1-\frac{1}{p(p-1)}\right)\nonumber\\
&=&\frac{1}{2^f}\sum_{l\geq0}\prod_{i=1}^{r}\frac{q_i+1}{q_i^{2e_i+1}}
\prod_{\scr p\ne q_i \atop \scr p:\mbox{\scrs odd}}
\left(1-\frac{1}{p(p-1)}\right).
\label{eq:5-4}
\end{eqnarray}
And, instead of $k=2^f(4l+1)$, we consider the similar sum for
$$k'=2^f(4l+3)=2^fq_1^{e_1}\cdots q_r^{e_r},$$
then we have
$$I^{(3)}(f)
=\frac{1}{2^f}\sum_{l\geq0}\prod_{i=1}^{r}\frac{q_i+1}{q_i^{2e_i+1}}
\prod_{\scr p\ne q_i \atop \scr p:\mbox{\scrs odd}}
\left(1-\frac{1}{p(p-1)}\right).$$
Thus
\begin{eqnarray*}
I^{(1)}(f)+I^{(3)}(f)&=&\frac{1}{2^f}
\left(\sum_{\scr 4l+1 \atop \scr l\geq0}+\sum_{\scr 4l+3 \atop \scr l\geq0}
\right)\prod_{i=1}^{r}\frac{q_i+1}{q_i^{2e_i+1}}
\prod_{\scr p\ne q_i \atop \scr p:\mbox{\scrs odd}}
\left(1-\frac{1}{p(p-1)}\right)\\
&=&\frac{1}{2^f}\prod_{p:\mbox{\scrs odd}}\left(1-\frac{1}{p(p-1)}\right)
\sum_{k:\mbox{\scrs odd}}\prod_{i=1}^{r}\left(\frac{q_i+1}{q_i^{2e_i+1}}
\cdot\frac{q_i(q_i-1)}{q_i^2-q_i-1}\right)\\
&=&\frac{1}{2^f}\prod_{p:\mbox{\scrs odd}}\left(1-\frac{1}{p(p-1)}\right)
\prod_{p:\mbox{\scrs odd}}\left\{1+\frac{(p+1)(p-1)}{p^2-p-1}
\sum_{j=1}^\infty \frac{1}{p^{2j}}\right\}\\
&=&\frac{1}{2^f},
\end{eqnarray*}
this proves Lemma \ref{lem:5-1}. \qed

\medskip
\noindent{\bf Proof of Lemma \ref{lem:5-2}.}

\medskip
\noindent Since $s$ is square free, $d$ is independent from the condition $s|\lcm{\ul{n},\ul{k}d}$ and (\ref{eq:5-2}) holds again. 

Now we put 
$$s=p_1\cdots p_t \quad \mbox{with $p_i$: odd prime}$$
and $(s,k_0)=z$. Then the condition ``$s|\lcm{\ul{n},\ul{k}d}$'' is equivalent to the condition ``$(s/z)|n$''. And in the same way to (\ref{eq:5-3}), we have
\begin{equation}
J^{(1)}(f,s)=\frac{1}{2^f}\sum_{l\geq0}\prod_{i=1}^r
\frac{q_i+1}{q_i^{2e_i}(q_i-1)}
\psum_{\scr n:\mbox{\scrs odd}\atop \scr \frac{s}{z}|n}\frac{\mu(n)}{n}
\prod_{\scr p|n \atop \scr p\ne q_i}\frac{1}{p-1}.
\label{eq:5-5}
\end{equation}
Here we introduce the sets, for $y|s$, 
\begin{eqnarray*}
{\bf N}^{(1)}_y:&=&\bigl\{k\in{\bf N}\ ;\ k\equiv1\pmod4,\ (k,s)=y \bigr\},\\
{\bf N}^{(3)}_y:&=&\bigl\{k\in{\bf N}\ ;\ k\equiv3\pmod4,\ (k,s)=y \bigr\}
\end{eqnarray*}
and 
$${\bf N}_y:={\bf N}^{(1)}_y\cup{\bf N}^{(3)}_y
=\bigl\{k:\mbox{odd}\ ;\ (k,s)=y \bigr\}.$$
Then
$$\bigcup_{y|s}{\bf N}^{(1)}_y=\bigl\{k\in{\bf N}\ ;\ k\equiv1\pmod4 \bigr\},$$
and this is a disjoint union. Now we calculate the partial sum of (\ref{eq:5-5}) over such a $k\in {\bf N}^{(1)}_y$. 

First we consider the sum, for a fixed $y|s$, 
\begin{equation}
S_y^{(1)}=\psum_{\scr n:\mbox{\scrs odd}\atop\scr \frac{s}{y}|n}
\frac{\mu(n)}{n}\prod_{\scr p|n \atop\scr p\ne q_i}\frac{1}{p-1}.
\label{eq:5-6}
\end{equation}
Let $n\in{\bf N}$ be odd and square free, we put
$$n=\frac{s}{y}\cdot m,$$
then
\begin{eqnarray*}
S_y^{(1)}&=&\psum_{\scr m:\mbox{\scrs odd}\atop\scr (m,\frac{s}{y})=1}
\frac{\mu\left(m\cdot\Frac{s}{y}\right)}{m\cdot\Frac{s}{y}}
\prod_{\scr p|m\cdot\frac{s}{y} \atop\scr p\ne q_i}\frac{1}{p-1}\\
&=&\frac{\mu(s)}{s}\cdot\frac{y}{\mu(y)}
\psum_{\scr m:\mbox{\scrs odd}\atop\scr (m,\frac{s}{y})=1}\frac{\mu(m)}{m}
\prod_{p|\frac{s}{y}}\frac{1}{p-1}
\prod_{\scr p|m \atop\scr p\ne q_i}\frac{1}{p-1}\\
&=&\left(\frac{\mu(s)}{s\varphi(s)}\right)
\left(\frac{\mu(y)}{y\varphi(y)}\right)^{-1}
\psum_{\scr m:\mbox{\scrs odd}\atop\scr (m,\frac{s}{y})=1}\frac{\mu(m)}{m}
\prod_{\scr p|m \atop\scr p\ne q_i}\frac{1}{p-1}\\
&=&\left(\frac{\mu(s)}{s\varphi(s)}\right)
\left(\frac{\mu(y)}{y\varphi(y)}\right)^{-1}
\prod_{\scr p\nmid \frac{s}{y} \atop{\scr p\ne q_i \atop\scr p:\mbox{\scrs odd}}}
\left(1+\frac{-1}{p}\cdot\frac{1}{p-1}\right)
\prod_{\scr p\nmid \frac{s}{y} \atop\scr p=q_i}\left(1+\frac{-1}{p}\right)\\
&=&\left(\frac{\mu(s)}{s\varphi(s)}\right)
\left(\frac{\mu(y)}{y\varphi(y)}\right)^{-1}\frac{\varphi(\ul{k_0})}{\ul{k_0}}
\prod_{\scr p\nmid \frac{s}{y} \atop{\scr p\ne q_i \atop\scr p:\mbox{\scrs odd}}}
\left(1-\frac{1}{p(p-1)}\right),
\end{eqnarray*}
and 
\begin{equation}
\sum_{\scr l\geq0 \atop\scr 4l+1\in{\bf N}_y^{(1)}}
\prod_{i=1}^r\frac{q_i+1}{q_i^{2e_i}(q_i-1)}S_y^{(1)}
=\left(\frac{\mu(s)}{s\varphi(s)}\right)
\left(\frac{\mu(y)}{y\varphi(y)}\right)^{-1}
\sum_{\scr l\geq0 \atop\scr 4l+1\in{\bf N}_y^{(1)}}
\prod_{i=1}^r\frac{q_i+1}{q_i^{2e_i+1}}
\prod_{\scr p\nmid \frac{s}{y} \atop{\scr p\ne q_i \atop\scr p:\mbox{\scrs odd}}}
\left(1-\frac{1}{p(p-1)}\right).
\label{eq:5-7}
\end{equation}
We can prove the similar formula for $\sum_{l\geq0,\ 4l+3\in{\bf N}_y^{(3)}}$ and now we consider the sum
$$T_y:=\left(\sum_{\scr l\geq0 \atop\scr 4l+1\in{\bf N}_y^{(1)}}
+\sum_{\scr l\geq0 \atop\scr 4l+3\in{\bf N}_y^{(3)}}\right)
\prod_{i=1}^r\frac{q_i+1}{q_i^{2e_i+1}}
\prod_{\scr p\nmid \frac{s}{y} \atop{\scr p\ne q_i \atop\scr p:\mbox{\scrs odd}}}
\left(1-\frac{1}{p(p-1)}\right).$$
Taking into account the following equivalence relation
$$k\in{\bf N}_y \quad\Leftrightarrow\quad \left\{\begin{array}{l}
\forall p|y \Rightarrow p|k\\
p|s \mbox{ and }p\nmid y \Rightarrow p\nmid k,
\end{array}\right.$$
we have
\begin{eqnarray*}
T_y&=&\prod_{p\nmid\frac{s}{y}}\left(1-\frac{1}{p(p-1)}\right)
\sum_{k\in{\bf N}_y}\prod_{i=1}^r\frac{q_i+1}{q_i^{2e_i+1}}\cdot
\frac{q_i(q_i-1)}{q_i^2-q_i-1}\\
&=&\prod_{p\nmid\frac{s}{y}}\left(1-\frac{1}{p(p-1)}\right)
\prod_{p|y}\left\{\frac{(p+1)(p-1)}{p^2-p-1}
\sum_{j=1}^\infty\frac{1}{p^{2j}}\right\}
\prod_{\scr p|s \atop p\nmid y}1\\
&\cdot& \prod_{p\nmid s}\left\{1+\frac{(p+1)(p-1)}{p^2-p-1}
\sum_{j=1}^\infty\frac{1}{p^{2j}}\right\}\\
&=&\prod_{p\nmid\frac{s}{y}}\left(1-\frac{1}{p(p-1)}\right)
\prod_{p|y}\frac{1}{p^2-p-1}\prod_{p\nmid s}\left(1+\frac{1}{p^2-p-1}\right).
\end{eqnarray*}
Therefore
\begin{eqnarray*}
T_y&=&\left(\frac{\mu(s)}{s\varphi(s)}\right)
\left(\frac{\mu(y)}{y\varphi(y)}\right)^{-1}\prod_{p|y}\frac{1}{p(p-1)}\\
&=&\frac{\mu(s)}{s\varphi(s)}\cdot \mu(y).
\end{eqnarray*}
Consequently, we have
\begin{eqnarray*}
J^{(1)}(f,s)+J^{(3)}(f,s)&=&\frac{1}{2^f}\sum_{y|s}T_y=\frac{1}{2^f}
\frac{\mu(s)}{s\varphi(s)}\sum_{y|s}\mu(y)\\
&=&\left\{\begin{array}{ll}
0, & \mbox{if }s>1,\\[0.3cm]
\Frac{1}{2^f}, & \mbox{if }s=1,
\end{array}\right.
\end{eqnarray*}
this proves Lemma \ref{lem:5-2}. \qed

\medskip
\noindent{\bf Proof of Lemma \ref{lem:5-3}.}

\medskip
\noindent We use the same notations with the proof of Lemma \ref{lem:4-2}. Lemma \ref{lem:5-1} and Lemma \ref{lem:5-2} show that
\begin{equation}
J^{(1)}(f,s)+J^{(3)}(f,s)=\left\{\begin{array}{ll}
\Frac{1}{2^f},& \mbox{if }s=1,\\[0.3cm]
0, & \mbox{if }s>1.
\end{array}\right.
\label{eq:5-8}
\end{equation}
Now we calculate the sum $J^{(1)}(f,s)-J^{(3)}(f,s)$, and we start from the formula (\ref{eq:5-7}). 

Here we construct the completely multiplicative function $\gamma$ by
$$\gamma(p)=\left\{\begin{array}{ll}
+1, & \mbox{if }p\equiv1\pmod4,\\
-1, & \mbox{if }p\equiv3\pmod4,\\
0, & \mbox{if }p=2,
\end{array}\right.$$
then, for an odd integer $k$, 
$$k\equiv3\pmod4 \quad\Leftrightarrow\quad \gamma(k)=-1.$$
Making use of this relation, we have
\begin{eqnarray*}
J^{(1)}(f,s)-J^{(3)}(f,s)&=&\frac{1}{2^f}\sum_{y|s}\frac{\mu(s)}{s\varphi(s)}
\left(\frac{\mu(y)}{y\varphi(y)}\right)^{-1}
\prod_{\scr p\nmid \frac{s}{y} \atop\scr p:\mbox{\scrs odd}}
\left(1-\frac{1}{p(p-1)}\right)\\
& &\cdot\sum_{k\in{\bf N}_y}\prod_{i=1}^r
\gamma(q_i^{e_i})\frac{q_i+1}{q_i^{2e_i+1}}\cdot\frac{q_i(q_i-1)}{q_i^2-q_i-1}.
\end{eqnarray*}
As for the sum $\sum_{k\in{\bf N}_y}$ in the right hand side, we have
\begin{eqnarray*}
\sum_{k\in{\bf N}_y}& & \!\!\!\!\!\!\!\!\!\!\!\prod_{i=1}^r\gamma(q_i^{e_i})
\frac{(q_i-1)(q_i+1)}{q_i^2-q_i-1}\cdot\frac{1}{q_i^{2e_i}}\\
&=&\prod_{\scr p|y \atop\scr \gamma(p)=1}
\left(\frac{(p-1)(p+1)}{p^2-p-1}\sum_{j=1}^\infty\frac{1}{p^{2j}}\right)
\prod_{\scr p|y \atop\scr \gamma(p)=-1}
\left(\frac{(p-1)(p+1)}{p^2-p-1}\sum_{j=1}^\infty\frac{(-1)^j}{p^{2j}}\right)
\prod_{\scr p|s \atop\scr p\nmid y}1\\
&\cdot&\prod_{\scr p\nmid s \atop\scr \gamma(p)=1}
\left\{1+\frac{(p+1)(p-1)}{p^2-p-1}\sum_{j=1}^\infty\frac{1}{p^{2j}}\right\}
\prod_{\scr p\nmid s \atop\scr \gamma(p)=-1}
\left\{1+\frac{(p+1)(p-1)}{p^2-p-1}\sum_{j=1}^\infty\frac{(-1)^j}{p^{2j}}\right\}.
\end{eqnarray*}
Thus
\begin{eqnarray*}
J^{(1)}(f,s)-J^{(3)}(f,s)&=&\frac{1}{2^f}\sum_{y|s}\frac{\mu(s)}{s\varphi(s)}
\left(\frac{\mu(y)}{y\varphi(y)}\right)^{-1}\\
&\cdot&\prod_{\scr p|y \atop\scr \gamma(p)=1}\frac{1}{p(p-1)}
\prod_{\scr p|y \atop\scr \gamma(p)=-1}\frac{-(p+1)}{p(p^2+1)}
\prod_{\scr p\nmid s \atop\scr \gamma(p)=-1}\frac{p^3-p^2-p-1}{(p^2+1)(p-1)}\\
&=&\frac{1}{2^f}\cdot\frac{\mu(s)}{s\varphi(s)}C\sum_{y|s}
\prod_{\scr p|y \atop\scr \gamma(p)=1}(-1)
\prod_{\scr p|y \atop\scr \gamma(p)=-1}\frac{(p+1)(p-1)}{p^2+1}\\
&\cdot& \prod_{\scr p|s \atop\scr \gamma(p)=-1}\frac{(p^2+1)(p-1)}{p^3-p^2-p-1}\\
&=&\frac{C}{2^f}\cdot\frac{\mu(s)}{s\varphi(s)}
\prod_{\scr p|s \atop\scr \gamma(p)=-1}\frac{2p^2(p-1)}{p^3-p^2-p-1}
\prod_{\scr p|s \atop\scr \gamma(p)=1}(1+\mu(p))\\
&=&\frac{C}{2^f}\prod_{\scr p|s \atop\scr \gamma(p)=1}\frac{-(1+\mu(p))}{p(p-1)}
\prod_{\scr p|s \atop\scr \gamma(p)=-1}\frac{-2p}{p^3-p^2-p-1}.
\end{eqnarray*}
This proves
$$J^{(1)}(f,s)-J^{(3)}(f,s)=\left\{\begin{array}{ll}
0, & \mbox{if $s$ has a prime divisor $p$}\\
 & \mbox{with }p\equiv1\pmod4,\\[0.3cm]
\Frac{C}{2^f}\prod_{p|s}\Frac{-2p}{p^3-p^2-p-1}, & 
\mbox{if all prime divisors $p$ of $s$ }\\[-0.4cm]
 & \mbox{satisfy }p\equiv3\pmod4.
\end{array}\right.$$
Combining with (\ref{eq:5-8}), we get Lemma \ref{lem:5-3}. \qed

\medskip
\noindent{\bf Proof of Theorem \ref{th:main}.}

\medskip
\noindent{\bf (I)} When $a_1\equiv1\pmod4$ and $a_1>1$, with Corollary \ref{cor:3-1} and Proposition \ref{prop:4-1} (1), (\ref{eq:1-2}) and (\ref{eq:1-3}) turn into
\begin{eqnarray*}
\delta_1=\delta_3
&=&\sum_{f\geq1}\frac{1}{2^{f+1}}\left(I^{(1)}(f)+J^{(1)}(f,a_1)\right)\\
&+&\sum_{f\geq1}\frac{1}{2^{f+1}}\left(I^{(3)}(f)+J^{(3)}(f,a_1)\right)
\end{eqnarray*}
and Lemma \ref{lem:5-1} and Lemma \ref{lem:5-2} show that
$$\delta_1=\delta_3=\sum_{f\geq1}\frac{1}{2^{2f+1}}=\frac16.$$
Clearly, the same argument yields the same result for the case $a_1\equiv3\pmod4$. 

\medskip
\noindent{\bf (II)} When $a_1\equiv2\pmod4$ with $a_1=2a_1'$. 

\medskip
\noindent{\bf Case (i)} $a_1'=1$.

\medskip
\noindent By Corollary \ref{cor:3-1} and Proposition \ref{prop:4-1} (1) again, we have from (\ref{eq:1-2}) that
\begin{eqnarray*}
\delta_1&=&\sum_{f\geq3}\sum_{l\geq0}\frac{k_0}{\varphi(k_0)}
\sum_{\scr d|k_0 \atop\scr d:\mbox{\scrs odd}}\frac{\mu(d)}{d}
\psum_{n}\frac{\mu(n)}{\frac12 nk\cdot 2^{f+1}\varphi(\lcm{\ul{n},\ul{k}d})}\\
&+&\sum_{\scr f\geq3 \atop\scr f=1}\sum_{l\geq0}\frac{k'_0}{\varphi(k'_0)}
\sum_{\scr d|k'_0 \atop\scr d:\mbox{\scrs odd}}\frac{\mu(d)}{d}
\psum_{n}
\frac{\mu(n)}{\frac12 nk'\cdot 2^{f+1}\varphi(\lcm{\ul{n},\ul{k'}d})}\\
&=&\sum_{f\geq3}\frac{1}{2^f}\left(I^{(1)}(f)+I^{(3)}(f)\right)
+\frac12I^{(3)}(1)
\end{eqnarray*}
Then Lemma \ref{lem:5-1} and Lemma \ref{lem:5-3} give
\begin{eqnarray*}
\delta_1&=&\sum_{f\geq3}\frac{1}{2^{2f}}+\frac12\cdot\frac{1-C}4\\
&=&\frac{7}{48}-\frac{C}{8}.
\end{eqnarray*}
Similarly, 
$$\delta_3=\frac{7}{48}+\frac{C}{8}.$$

\medskip
\noindent{\bf Case (ii)} $a_1'\equiv1\pmod4$ with $a'_1>1$.

\medskip
\noindent Applying Corollary \ref{lem:3-1} and Proposition \ref{prop:4-1} (2) (a) to (\ref{eq:1-2}), we can prove that
\begin{eqnarray*}
\delta_1&=&\sum_{f\geq1}\frac{1}{2^{f+1}}\left(I^{(1)}(f)+J^{(1)}(f,a_1')
+I^{(3)}(f)+J^{(3)}(f,a_1')\right)\\
&-&\sum_{f\geq1}\frac{1}{2^{f+1}}\cdot 2
\left(J^{(1)}(f,a_1')+J^{(3)}(f,a_1')\right)\\
&+&\sum_{f\geq3}\frac{2}{2^{f+1}}\left(J^{(1)}(f,a_1')+J^{(3)}(f,a_1')\right)
+2\cdot\frac14 J^{(3)}(1,a_1')\\
&=&\frac16+\frac12 J^{(3)}(1,a_1'),
\end{eqnarray*}
and Lemma \ref{lem:5-3} gives the desired result. 

\medskip
\noindent{\bf Case (iii)} $a_1'\equiv3\pmod4$.

\medskip
\noindent The same argument gives
$$\delta_1=\frac16+\frac12 J^{(1)}(1,a_1')$$
and Lemma \ref{lem:5-3} gives the desired result. \qed

\end{document}